\magnification\magstep1
\overfullrule=0pt
\voffset=-0.3 in
\input amstex
\documentstyle{amsppt}
\catcode`\@=11 \def\logo@{} \catcode`\@=\active 
\rightheadtext{ Omori-Yau maximum principle } \leftheadtext {A.
Borb\'ely} \pageno 1 \TagsOnRight
\def\({\left(}
\def\){\right)}
\def\[{\left[}
\def\]{\right]}

\def\<{\langle }
\def\>{\rangle }

\def\Ricci{\hbox{Ricc}}
\topmatter

\title     A remark on the Omori-Yau maximum principle
\endtitle
    \author Albert Borb\'ely     \endauthor

    \address Kuwait University, Faculty of Science,
         Department of Mathematics,
         P.O. Box 5969, Safat 13060, Kuwait
    \endaddress

\email borbely\@sci.kuniv.edu.kw \endemail


    \abstract   A Riemannian manifold  $M$ is said to satisfy the Omori-Yau maximum
    principle if for any $C^2$ bounded
    function $g:M\to \Bbb R$ there is a sequence $x_n\in M$, such
    that $\lim _{n\to \infty }g(x_n)=\sup _M g$, $ \lim _{n\to \infty }|\nabla
    g(x_n)|=0$ and $\limsup _{n\to \infty }\Delta g(x_n)\leq 0$.
    It is shown that if the Ricci curvature does not approach $-\infty $ too
    fast the manifold satisfies the Omori-Yau maximum principle.
    This improves earlier necessary conditions.
    The given condition is quite optimal.

   \endabstract
\endtopmatter
\document

\head    0. Introduction
\endhead

\proclaim {Definition} A Riemannian manifold  $M$ is said to
satisfy the Omori-Yau maximum
    principle if for any $C^2$
    function $g:M\to \Bbb R$ which is bounded from above and for
    any $\epsilon >0$ there is a point $x_{\epsilon }\in M$, such
    that $|g(x_{\epsilon })-\sup _M g|<\epsilon $, $ |\nabla
    g(x_{\epsilon })|<\epsilon $ and $\Delta g(x_n)<\epsilon $.
\endproclaim

This principle has turned out to be very useful in differential
geometry and received considerable attention recently. A necessary
condition in terms of the Ricci curvature for a manifold to
satisfy this principle was first proved by Omori in [O] and later
generalized by Yau [Y]. It states that if the Ricci curvature is
bounded from below then the manifold satisfies the Omori-Yau
maximum principle.

This was improved upon by Ratto, Rigoli and Setti in [RRS, Theorem
2.3].

\proclaim {Theorem (Ratto-Rigoli-Setti)} Let $M^n$ be a complete
Riemannian manifold, $p\in M^n$ be a fixed point and $r(x)$ be the
distance function from $p$. Let us assume that away from the cut
locus of $p$ we have $$\Ricci (\nabla r,\nabla r)\geq -
(n-1)BG^2(r),$$ where $B>0$ is some constant and $G(t)$ has the
following properties:
$$ \eqalign {&(i)\quad G(0)=1,\qquad  \quad G'\geq 0 \cr  \noalign {\vskip 4pt}&(ii)\quad \int_0^{\infty
} {dt\over G(t)}=\infty \cr \noalign {\vskip 4pt}&(iii)\quad
{d^{2k+1}\over dt^{2k+1}}\sqrt G(0)=0\quad \hbox{for all}\quad
k\in \Bbb N\cr \noalign {\vskip 4pt}&(iv)\quad \limsup _{t\to
\infty }{t\sqrt {G(\sqrt t)}\over \sqrt {G(t)}}<\infty .}$$ Then
$M^n$ satisfies the Omori-Yau maximum principle.
\endproclaim

The goal of the present note is to improve the necessary condition
given in [RRS]. The actual statement is given as a Corollary.
Basically we remove the last two conditions on the function
$G(t)$, which turned out not to be essential.

Another interesting necessary condition, requiring the existence
of an exhaustion function with certain properties, was given by
Kim and Lee in [KL]. Interestingly there is an alternative proof
by Kim and Lee  of the Ratto-Rigoli-Setti result in [KL] which is
still using these extra conditions.

The proof uses the same method we used in an earlier paper [B].

\proclaim {Theorem}  Let $M^n$ be a complete Riemannian manifold,
$p\in M^n$ be a fixed point and $r(x)$ be the distance function
from $p$. Let us assume that $$\Delta r(x)\leq G(r(x))$$ for all
$x\in M^n$ where $r$ is smooth and $r(x)>1$, where $G(t)$ has the
following properties:
$$ G\geq 1,\quad G'\geq 0,\quad \hbox {and}\quad \int_0^{\infty
} {dt\over G(t)}=\infty .$$ Then $M^n$ satisfies the Omori-Yau
maximum principle.
\endproclaim

As a consequence we have the following.

\proclaim {Corollary}  Let $M^n$ be a complete Riemannian
manifold, $p\in M^n$ be a fixed point and $r(x)$ be the distance
function from $p$. Let us assume that away from the cut locus of
$p$ we have $$\Ricci (\nabla r,\nabla r)\geq - G^2(r),$$ where
$G(t)$ has the following properties:
$$ G\geq 1,\quad G'\geq 0,\quad \hbox {and}\quad \int_0^{\infty
} {dt\over G(t)}=\infty .$$ Then $M^n$ satisfies the Omori-Yau
maximum principle.
\endproclaim

The main condition on the function $G(t)$ in the Corollary and in
 the Ratto-Rigoli-Setti Theorem  is the same ($\int 1/G(t)=\infty $) but
 there are additional technical conditions imposed on the function
 $G(t)$ in the later Theorem.  In this respect Corollary can be considered as a refinement of
 the Ratto-Rigoli-Setti Theorem.

Let us mention that  this condition  is quite optimal. If
$\int_0^{\infty } 1/G(t)dt<\infty $, there are manifolds with
$\Delta r\leq G(r) $ for which the Omori-Yau maximum principle
does not apply. The details can be found in Section 3.

\head 1. Proof of the Theorem \endhead

\demo {Proof of the Theorem} Set $L=\sup g$ and let us assume that
$g<L$ at every point of $M$. Otherwise $g$ assumes its maximum at
some point and that point trivially satisfies the conditions of
the Definition for all $\epsilon >0$.

Define the function $F(t)$ as
$$F(t)=e^{\int _0^{t }{1\over G(s)}ds}.$$
Then clearly:\ \  $\displaystyle F\geq 1,\quad F\hbox {  is
strictly increasing and}\quad \lim _{t\to \infty }F(t)=\infty .$

For any $\epsilon <  \min \{1,L-\sup \{g(x):r(x)<1\}\}$ define the
function $h_{\lambda }:M\to \Bbb R$ as
$$ h_{\lambda }(x)=\lambda F(r(x))+L-\epsilon .$$

Since $F(r(x))\geq 1 $, for $\lambda
>\epsilon $ we have

$$h_{\lambda }(x)>L>g(x)\qquad \hbox { for all}\quad  x\in M.$$

Define $\lambda _0$ as
$$ \lambda _0=\inf \{\lambda :h_{\lambda }(x)>g(x)\quad \hbox {for
all}\quad x\in M\}.$$ Since $\sup g=L$ it is easy to see that
$\lambda _0>0$ and $h_{\lambda _0}(x)\geq g(x)$ for all $x\in M$.

We claim that there is a point $x_{\epsilon }\in M$ such that
$h_{\lambda _0}(x_{\epsilon })= g(x_{\epsilon })$.

This will follow from the  observation that if $h_{\lambda }(x)>
g(x)$ for all $x\in M$, then there is a $\lambda '<\lambda $ such
that $h_{\lambda '}(x)> g(x)$ for all $x\in M$. To show this we
argue as follows.

Let $r_0$ be large enough such that $h_{\lambda }(x)> L+1$ for
$r(x)>r_0.$ Since $\lim _{r\to \infty }F(r)=\infty $ such $r_0$
must exists. The set $\{x\in M: r(x)\leq r_{0 }\}$ is compact,
therefore $h_{\lambda }(x)> g(x)$ for all $x\in M$ implies that
there is a $\lambda '<\lambda $ such that $h_{\lambda '}(x)> g(x)$
for all $x\in M: r(x)\leq r_{0 }\}$. Choosing $\lambda '$
sufficiently close to $\lambda $ we can achieve that $h_{\lambda '
}(x)> L$ for $r(x)=r_0.$ Since $F$ is increasing we obtain that
$h_{\lambda ' }(x)> L$ for $r(x)\geq r_0.$ Combining this with the
previous remark we have $h_{\lambda ' }(x)> g(x)$ for all $x\in
M.$

Next, we have to show that $h_{\lambda _0}$ is smooth at
$x_{\epsilon }$. The argument is exactly the same as the argument
 in [B], but we include it at the and of this proof for the
convenience of the reader.

Once we established the smoothness of $h_{\lambda _0}$ at
$x_{\epsilon }$, the rest of the argument is straight forward.

From the definition of $F$ and from the fact that $G'\geq 0$ we
have
$$F'={F\over G}\qquad \hbox{and}\qquad F''={F'\over G}-{FG'\over
G^2}\leq {F\over G^2}.$$

From the fact that $g(x_{\epsilon })=\lambda _0F(r(x_{\epsilon
}))+L-\epsilon <L$ we conclude that
$$L-g(x_{\epsilon })\leq \epsilon ,\tag 1.1$$
moreover
$$\lambda _0F(r(x_{\epsilon }))<\epsilon \qquad \hbox {hence}  \qquad \lambda
_0<{\epsilon \over F(r(x_{\epsilon }))}<\epsilon .\tag 1.2$$

Since
$$ h_{\lambda _0}(x)\geq g(x),\quad \hbox {and }\quad h_{\lambda _0}(x_{\epsilon })=g(x_{\epsilon
}),$$ we have
$$\nabla g(x_{\epsilon })=\nabla h_{\lambda _0}(x_{\epsilon
})\quad \hbox {and}\quad \Delta h_{\lambda _0}(x_{\epsilon })\geq
\Delta g(x_{\epsilon }).$$

Taking into consideration (1.2), the definition of $F$, the fact
that $|\nabla r|=1$ and the assumption that $G(r)\geq 1$, the
first equality above yields
$$|\nabla g(x_{\epsilon })|=|\lambda _0F'(r(x_{\epsilon }))\nabla r(x_{\epsilon })|=
{\epsilon \over F(r)}\cdot {F(r)\over G(r)}<\epsilon .\tag 1.3$$

For the Laplace of $h_{\lambda _0}$  we have
$$\eqalign {\Delta g(x_{\epsilon })&\leq \Delta h_{\lambda _0}(x_{\epsilon })=
\lambda _0\bigg (F'(r(x_{\epsilon }))\Delta r(x_{\epsilon
})+F''(r(x_{\epsilon }))|\nabla r(x_{\epsilon })|^2\bigg )\leq \cr
&\leq {\epsilon \over F}\bigg ({F\over G}\Delta r+{F\over
G^2}\bigg )\leq 2\epsilon .}\tag 1.4$$

The inequalities (1.1), (1.3) and (1.4) show that the point
$x_{\epsilon }$ satisfies the conditions in the Definition.

 Finally, we have to show that $h_{\lambda _0}$ is smooth
at $x_{\epsilon }$. Since $ h_{\lambda }(x)=\lambda
F(r(x))+L-\epsilon $ it is enough to show that $r$ is smooth at
$x_{\epsilon }$. If not, then $x_{\epsilon }$ must be on the cut
locus of $p$. In this case we have two possibilities. Either there
are two distinct minimizing
 geodesic segments $\gamma _1,\gamma _2:[0,t_0]\to M$
joining $p$ to $x_{\epsilon }$, or there is a  geodesic segment
$\gamma :[0,t_0]\to M$  from $p$ to $x_{\epsilon }$ along which
$x_{\epsilon }$ is conjugate to $p$.

In both cases we have $$t_0=r(x_{\epsilon }).$$

Let us start with the first case. Let $w=\gamma _1'(t_0)$ and
$v=\gamma _2'(t_0)$. Since $\gamma _1$ and $\gamma _2$ are
distinct segments we have $w\ne v$. The functions $t\to r(\gamma
_i(t))$ are differentiable on $(0,t_0)$ (for $i=1,2$) and they
have a left-derivative at $t_0$.

From the fact that $h_{\lambda _0}\geq g$ and $h_{\lambda
_0}(x_{\epsilon })=g(x_{\epsilon })$ we have
$$\liminf _{s\to 0^+}{h_{\lambda _0}(\gamma _2(t_0+s))-h_{\lambda _0}(\gamma _2(t_0))\over
s}\geq D_vg(x_{\epsilon }),$$ where $D_vg(x_{\epsilon })$ denotes
the directional derivative of $g$ at the point $x_{\epsilon }$ in
the direction of $v$. Moreover since $g$ is smooth and $h_{\lambda
_0}$ has a directional derivative at $x_{\epsilon }$ in the
direction of $-v$, we also have
$$-\lambda _0F'(r(x_{\epsilon }))=D_{-v}h_{\lambda _0}(x_{\epsilon })\geq D_{-v}g(x_{\epsilon
})=-D_vg(x_{\epsilon }).$$ This yields
$$D_vg(x_{\epsilon})\geq \lambda _0F'(r(x_{\epsilon })).\tag 1.5$$

Combining this with the above inequality we obtain
$$\liminf _{s\to 0^+}{h_{\lambda _0}(\gamma _2(t_0+s))-h_{\lambda _0}(\gamma _2(t_0))\over
s}\geq \lambda _0F'(r(x_{\epsilon })).$$ Taking into account the
special form of $h_{\lambda _0}$ we have
$$\liminf _{s\to 0^+}{r(\gamma _2(t_0+s))-r(\gamma _2(t_0))\over
s}\geq 1.\tag 1.6$$

This will lead to a contradiction. Since $v\ne w$, there is a
$0<c<1$ depending only on the angle of $v$ and $w$, such that
$$r(\gamma _2(t_0+s))<t_0+cs,\tag 1.7 $$ for a small enough $s>0$.

One can see this by connecting the point $\gamma _1(t_0-s)$ to
$\gamma _2(t_0+s)$ by a geodesic segment. Since $\gamma _1$ and
$\gamma _2$ are different there is a $0<c_1<1$ such that for a
small enough $s>0$ we have $\hbox{dist} (\gamma _1(t_0-s),\gamma
_2(t_0+s))<c_12s$ and this implies (1.7). Since $r(x_{\epsilon
})=r(\gamma _2(t_0))=t_0$ it is easy to see that (1.6) and (1.7)
are in direct contradiction.

We now turn our attention to the second case. Since $\gamma $ is
distance minimizing between $p$ and $x_{\epsilon }$ the distance
function $r$ is smooth at $\gamma (t)$ for $0<t<t_0$. Set
$m(t)=\Delta r(\gamma (t))$. Then $m(t)$ is also smooth on the
interval $(0,t_0)$ and since $\gamma (t_0)$ is conjugate to
$p=\gamma (0)$ along $\gamma $ we have
$$\lim _{t\to t_0^-}m(t)=-\infty .\tag 1.8$$

Since $\lambda _0 >0$, from (1.5) we conclude that
$D_vg(x_{\epsilon })>0$, that is $\nabla g(x_{\epsilon })\ne 0$.
This implies that the level surface $H=\{x\in M:
g(x)=g(x_{\epsilon })\}$ is a smooth hyper-surface near
$x_{\epsilon }$. Denote by $H_s$ the surface parallel to $H$ and
passing through the point $\gamma (t_0-s)$ for some $s>0$. Again,
since $H$ is smooth near $x_{\epsilon }$ the surface $H_s$ will
also be smooth near $\gamma (t_0-s)$ for a small enough  $s>0$.

It is now clear from (1.8) that for some  small $s>0$ we have
$$m(t_0-s)<\hbox {trace of the 2nd fundamental form of\ } H_s \hbox {\ at\ } \gamma
(t_0-s),$$ where the second fundamental form of $H_s$ at $\gamma
(t_0-s)$ is taken in the direction of $\gamma '(t_0-s)$.

Taking into account that  $m(t_0-s)$ is the trace of the 2nd
fundamental form of the geodesic ball $B_p(t_0-s)$ around $p$ at
the point $\gamma (t_0-s)$ (with respect to the same normal vector
$\gamma '(t_0-s)$) we conclude that there has to be a point $q_s
\in H_s$, sufficiently close to $\gamma (t_0-s)$, that lies inside
$B_p(t_0-s)$. This means that
$$r(q_s)<t_0-s.$$

Since $H_s$ is parallel to $H$ we have a point on $q\in F$ such
that $\hbox {dist}(q_s,q)=s$. Combining this with the above
inequality we have
$$ r(q)<t_0=r(x_{\epsilon }).$$
Taking into account that $F$ is strictly increasing we obtain
$$h_{\lambda _0}(q)=\lambda _0F(r(q))+L-\epsilon <\lambda
_0F(r(x_{\epsilon }))+L-\epsilon = h_{\lambda _0}(x_{\epsilon
})=g(x_{\epsilon })=g(q).$$ This leads to a contradiction since
$h_{\lambda _0}\geq g$ on $M$.

\enddemo

\head 2. Proof of the Corollary \endhead

Let $q\in M$ be a point away from the cut locus of $p$ and $\gamma
$ be a geodesic segment parameterized by arc length connecting $p$
to $q$. Set $m(t)=\Delta r(\gamma (t))$ and $R(t)=\Ricci (\gamma
'(t),\gamma '(t))$. Then it is well known that $m(t)$ satisfies
the Riccati inequality along $\gamma $. Taking into consideration
the condition on the Ricci curvature we have
$$m'(t)\leq -R(t)-{m^2(t)\over n-1}\leq G^2(t)-{m^2(t)\over
n-1}.$$

This implies that  $m$ is decreasing as long as $m>\sqrt {n-1}G$
and a simple argument shows that
$$m(t)<(\sqrt {n-1}+1)G,$$
for all $t>t_0$, where $t_0$ is a sufficiently large constant,
independent of  $G$.

This yields $$\Delta r< (\sqrt {n-1}+1)G\quad \hbox {if} \quad
r>t_0,
$$ for points that are not on the cut locus of $p$. Since $(\sqrt
{n-1}+1)G$ satisfies the conditions in the Theorem the proof of
the corollary is complete.

\head 3. An example \endhead

In this section we sketch an example, that shows that the
condition in the Theorem is quite optimal. Let  $M^n$ be a
Hadamard manifold that is rotationally symmetric around $p\in
M^n$.

 Let $r$ be the distance function from $p$ and
assume that $\Delta r(x)>G(r)$ for all $x\in M^n$, where $G$
satisfies the conditions:$$ G\geq 1,\quad G'\geq 0,\quad \hbox
{and}\quad \int_0^{\infty } {dt\over G(t)}<\infty .$$

Then there is a bounded function $h:M\to \Bbb R$ which shows that
the manifold $M^n$ does not satisfy the Omori-Yau maximum
principle. To construct $h$ we need the following lemma.

\proclaim {Lemma} Let $G:[0,\infty )\to \Bbb R$ be a function
satisfying the conditions: $$ G\geq 1,\quad G'\geq 0,\quad \hbox
{and}\quad \int_0^{\infty } {dt\over G(t)}<\infty .$$ Then there
is a function $H:[0,\infty )\to \Bbb R$ such that $$ H\geq
1/2,\quad H'\geq 0,\quad 2H\leq G,\quad H'\leq H^2 \quad \hbox
{and}\quad \int_0^{\infty } {dt\over H(t)}<\infty .$$
\endproclaim

First we  construct the function $h:M^n\to \Bbb R$ and give the
proof of the Lemma later.

Let
$$  h(x)=\int _0^{r(x)}{dt\over H(t)}.$$

The last condition on $H$ in the Lemma implies that $h$ is bounded
from above. A simple computation shows that
$$\Delta h={\Delta r\over H}-{H'\over H^2}|\nabla r|^2.$$
Since $\Delta r>G(r)\geq 2H(r)$, $|\nabla r|=1$ and $H'\leq H$ we
have
$$\Delta h>2-{1}=1.$$
This clearly shows that the manifold $M^n$ does not satisfy the
Omori-Yau maximum principle.

All that remains is to prove the Lemma.

\demo {Proof of Lemma} Let $A\subset (0,\infty )$ be defined as
$$A=\{ t>0: {G'(t)\over 2}>\bigg({G(t)\over 2}\bigg)^2\}.$$
It is an open set therefore
$$A=\cup I_n,$$
where $I_n=(t_n,s_n)$ are disjoint open intervals.

This is the set where $G/2$ grows too fast. We obtain $H$ by
modifying $G/2$ on a slightly larger set so that it will never
grow too fast, that is $H'\leq H^2$.

For a given $n$ define the function $k_n(t)$ to be
$$k_n(t)={1\over a_n-t},$$
where $a_n$ is chosen such that $k_n(t_n)=G(t_n)/2$. Then we have
$$k_n(t_n)={G(t_n)\over 2},\quad k_n'(t)=k_n^2(t)\quad \hbox {and}\quad
{G'(t)\over 2}> \bigg({G(t)\over 2}\bigg)^2\quad \hbox {for}\quad
t\in (t_n,\min \{s_n,a_n\}).$$ This implies that
$$k_n(t)<{G(t)\over 2}\quad  \hbox
{for}\quad t\in (t_n,\min \{s_n,a_n\}).$$

Let $v_n>t_n $ be the first point where $k_n(v_n)=G(v_n)/2$. Such
point must exists since $\lim _{t\to a_n}k_n(t)=\infty .$
Therefore we have $t_n<s_n<v_n<a_n$ and as a result
$J_n=(t_n,v_n)\supset I_n$.

The intervals $I_n$ are all disjoint but $J_n$ are not necessarily
disjoint intervals. However if $J_n\cap J_m\ne \emptyset $, then
either $J_n\subset J_m$ or $J_m\subset J_n$. This follows simply
from the way the intervals $J_n$ were constructed and from the
fact that the graphs of the functions $1/(a-t)$, $t<a$ and
$1/(b-t)$, $t<b$ are translates of each other.

Therefore we can select a pairwise disjoint family of intervals
$J_{n_l}$ such that $B=\cup J_n=\cup J_{n_l}$. To simplify the
notation without loss of generality we can assume that the
intervals $J_n$ are already pairwise disjoint.

We can now define the function $H(t)$ as follows
$$ H(t)= \left \{ \quad \eqalign { &{G(t)\over 2}\cr
&{1\over a_n-t}\cr } \qquad \eqalign { &\hbox {if}\quad  t\notin
B=\cup J_n\cr \noalign {\vskip 10pt} &\hbox{if}\quad t\in J_n.\cr
}\right.$$

It is clear from the construction that $H$ satisfies the first
four properties in the Lemma. It remains to show that it will
satisfy the remaining property
$$\int _0^{\infty }{dt\over H(t)}<\infty.\tag 3.1 $$
We can write
$$\int _0^{\infty }{dt\over H(t)}=\int _B{dt\over H(t)}+\int _{\Bbb R^+-B}{dt\over H(t)}.$$

The second integral is clearly finite since
$$\int _{\Bbb R^+-B}{dt\over H(t)}=\int _{\Bbb R^+-B}{2dt\over
G(t)}<\infty .$$

The first integral can be computed as follows
$$\int _B{dt\over H(t)}=\int _{\cup J_n}{dt\over H(t)}=\sum
_{n=1}^{\infty }\int _{t_n}^{v_n}{ a_n-t}\,\, dt={1\over 2}\sum
_{n=1}^{\infty }(a_n-t_n)^2-(a_n-v_n)^2.$$

From the construction of the intervals $J_n$ and the function $H$
one obtains that
$$a_n-t_n={2\over G(t_n)}\qquad \hbox {and}\qquad a_n-v_n={2\over
G(v_n)}.$$

To show that the infinite sum above is finite it is enough to show
that any partial-sum is bounded by a fixed constant. For this
reason consider the sum
$$ \sum
_{n=1}^{m}(a_n-t_n)^2-(a_n-v_n)^2=\sum _{n=1}^{m}\bigg( {2\over
G(t_n)}\bigg)^2-\bigg( {2\over G(v_n)}\bigg)^2 .$$ By rearranging
the terms if necessary, without loss of generality we can assume
that
$$t_1<v_1<t_2<v_2<....<t_n<v_n<t_{n+1}<v_{n+1}<...<t_m<v_m.$$
Taking into consideration that $G(t)$ is an increasing function,
we obtain that
$$ \sum
_{n=1}^{m}(a_n-t_n)^2-(a_n-v_n)^2=\sum _{n=1}^{m}\bigg( {2\over
G(t_n)}\bigg)^2-\bigg( {2\over G(v_n)}\bigg)^2 <\bigg( {2\over
G(t_1)}\bigg)^2.$$ This shows that the above sum is finite, which
in turn proves (3.1). This completes the proof of the lemma.

\enddemo


\Refs \widestnumber\key{MNP}

\ref \key B \by A. Borb\'ely\paper  Immersion of manifolds with
unbounded image and a modified maximum principle of Yau \jour
Bull. Australian Math. Soc. \vol 78 \yr 2008 \pages 285-291
\endref

\ref \key KL \by K. Kim, H. Lee \paper On the Omori-Yau almost
maximum principle \jour J. Math. Anal. Appl. \vol 335 \yr 2007
\pages 332-340
\endref

\ref \key O \by H. Omori \paper Isometric immersions of Riemannian
manifolds \jour J. Math. Soc. Jpn. \vol 19 \yr 1967 \pages 205-214
\endref

\ref \key RRS \by A. Ratti, M. Rigoli, A.G.Setti \paper On the
Omori-Yau maximum principle and its application to differential
equations and geometry \jour J. Funct. Anal. \vol 134 \yr 1995
\pages 486-510
\endref

\ref \key Y \by S.-T. Yau \paper Harmonic functions on complete
Riemannian Manifolds \jour Communications on pure and applied
Mathematics \vol 28 \yr 1975 \pages 201-228 \endref







\endRefs
\enddocument